\documentclass[11pt]{article}
 \usepackage[T1]{fontenc}
\usepackage{amsthm}
\usepackage{amssymb}
\usepackage{amsmath}
\usepackage{comment}
\usepackage{thm-restate}
\usepackage{latexsym}
\usepackage{graphicx,enumerate,color}
\usepackage{url} 
\usepackage[affil-it]{authblk}
\usepackage{hyperref}
\usepackage[noabbrev,capitalise]{cleveref}

\usepackage{color}
\usepackage[normalem]{ulem}

\usepackage[margin=1in]{geometry}

\newcommand{\eps}{\varepsilon}

\theoremstyle{plain}
\newtheorem{thm}{Theorem}[section]
\newtheorem{lem}[thm]{Lemma}

\newtheorem{conj}[thm]{Conjecture}

\newtheorem{prob}[thm]{Problem}

{\noindent \emph{Proof.} {}{#1}{}}{\hfill
	$\Diamond$\vspace{1em}}

\theoremstyle{plain} 
\newcommand{\thistheoremname}{}
\newtheorem{genericthm}[section]{\thistheoremname}

\theoremstyle{definition}

\newcounter{counter}
\newcommand{\less}{\setminus}

 \def\es{\emptyset}
 \def\dfn#1{{\sl #1}}
 
\begin{document}
\date{}
 
\title{On the size of $(K_t, K_{1,k})$-co-critical graphs}
\author{Hunter Davenport$^{a,}$\thanks{Supported by The Honors in the Major Program for undergraduate students at the University of Central Florida.}\hskip 1cm   Zi-Xia Song$^{a,}$\thanks{Supported by  NSF award  DMS-1854903.  E-mail address:    Zixia.Song@ucf.edu}\hskip 1cm  Fan Yang$^{b}$}

 \affil{ 
  { \small $^a${Department  of Mathematics, University of Central Florida, Orlando, FL 32816, USA}} 
   \\{ \small $^b${Department of Mathematics, Shandong University, Jinan, Shandong, 250100, China.}}  
     }

\maketitle

\begin{abstract}
Given   graphs $G, H_1, H_2$, we write \emph{$G \rightarrow ({H}_1, H_2)$} if  every $\{$red, blue$\}$-coloring of the edges of $G$ contains  a  red copy of $H_1$ or a blue copy of $H_2$. A non-complete graph $G$ is \emph{$(H_1, H_2)$-co-critical} if $G  \nrightarrow ({H}_1, H_2)$,   but  $G+e\rightarrow ({H}_1, H_2)$ for every edge $e$ in $\overline{G}$.  Motivated by a conjecture of   Hanson and Toft from 1987,      we study the minimum number of edges over all  $(K_t, K_{1,k})$-co-critical graphs on $n$ vertices.   
We  prove   that for all   $t\ge3$ and      $k\ge 3$,    there exists a constant $\ell(t, k)$ such that, for all      $n \ge (t-1)k+1$, if $G$  is  a $(K_t, K_{1,k})$-co-critical graph  on $n$ vertices, then    $$   e(G)\ge   \left(2t-4+\frac{k-1}{2}\right)n-\ell(t, k).$$ 
Furthermore, this linear bound is asymptotically best possible  when $t\in\{3, 4,5\}$ and   all $k\ge3$ and $n\ge (2t-2)k+1$. It seems non-trivial to construct extremal $(K_t, K_{1,k})$-co-critical graphs for $t\ge6$.  We also  obtain  the sharp bound for the size of $(K_3, K_{1,3})$-co-critical graphs  on $n\ge13$ vertices by showing that all such graphs have at least $3n-4$ edges.  \end{abstract}


\baselineskip=18pt

\section{Introduction}

 All graphs considered in this paper are finite, and without loops or multiple edges. For a graph $G$, we will use $V(G)$ to denote the vertex set, $E(G)$ the edge set, $|G|$ the order $|V(G)|$ of $G$, $e(G)$ the size $|E(G)|$ of $G$,   $N_G(x)$   the neighborhood of vertex $x$ in $G$, $\delta(G)$ the minimum degree, $\Delta(G)$ the maximum degree,  $\alpha(G)$ the independence number, and $\overline{G}$ the complement of $G$.
If  $A, B \subseteq V(G)$ are disjoint, we say that $A$ is \emph{complete to} $B$ if  every vertex in $ A$ is adjacent to every vertex in   $B$; and $A$ is  \emph{anti-complete to} $B$ if   no vertex in $A$ is adjacent to a vertex in $B$.
The subgraph of $G$ induced by $A$, denoted $G[A]$, is the graph with vertex set $A$ and edge set $\{xy \in E(G): x, y \in A\}$. We denote by $B \less A$ the set $B - A$, $e_G(A, B)$ the number of edges between $A$ and $B$ in $G$, and $G \less A$ the subgraph of $G$ induced on $V(G) \less A$, respectively.
If $A = \{a\}$, we simply write $B \less a$, $e_G(a, B)$, and $G \less a$, respectively.  For  any  edge  $e$ in $ \overline{G} $, we use $G+e$ to denote the graph obtained from $G$ by adding the new edge $e$. 
The {\emph{join}} $G+H$ (resp.~{\emph{union}} $G\cup H$) of two 
vertex-disjoint graphs
$G$ and $H$ is the graph having vertex set $V(G)\cup V(H)$  and edge set $E(G)
\cup E(H)\cup \{xy:   x\in V(G),  y\in V(H)\}$ (resp. $E(G)\cup E(H)$).
Given two isomorphic graphs $G$ and $H$, we may (with a slight but common abuse of notation) write $G = H$.   For an integer $t\ge1$ and a graph $H$, we define $tH$ to be the union of $t$ disjoint copies of $H$.  We use $K_n$, $K_{1,{n-1}}$,   $P_n$ and $T_n$ to denote the complete graph,  star,     path and a tree on  $n$ vertices, respectively.  For any positive integer $r$, we write  $[r]$ for the set $\{1,2, \ldots, r\}$. We use the convention   ``$A:=$'' to mean that $A$ is defined to be the right-hand side of the relation.  \medskip

Given an integer $k \ge 1$ and   graphs $G$, ${H}_1, \dots, {H}_k$, we write \dfn{$G \rightarrow ({H}_1, \dots, {H}_k)$} if every $k$-coloring of $E(G)$ contains a monochromatic  ${H}_i$ in color $i$ for some $i\in [k]$.
The classical \dfn{Ramsey number}  $r({H}_1, \dots, {H}_k)$  is the minimum positive integer $n$ such that $K_n \rightarrow ({H}_1, \dots, {H}_k)$.  Following Ne$\check{s}$et$\check{r}$il~\cite{Nesetril1986}, and Galluccio, Simonovits and Simonyi~\cite{Galluccio1992}, a non-complete graph $G$   is \emph{$(H_1, \ldots,  H_k)$-co-critical} if $G  \nrightarrow ({H}_1,  \ldots, {H}_k)$,   but  $G+e\rightarrow ({H}_1,  \ldots, {H}_k)$ for every edge $e$ in $\overline{G}$.    We simply  say  $G$ is $(H; k)$-co-critical  when  $H_1=\cdots=H_k=H$. \medskip

 The notation of co-critical graphs was initiated by Ne$\check{s}$et$\check{r}$il~\cite{Nesetril1986} in 1986. It is simple to check that   $K_6^-$ is $(K_3, K_3)$-co-critical, where   $K_6^-$ denotes  the graph obtained from $K_6$ by deleting exactly one edge, 
 and  every   $(H_1, \ldots,  H_k)$-co-critical graph has at least $r({H}_1, \dots, {H}_k)$   vertices.     Galluccio, Simonovits and Simonyi~\cite{Galluccio1992} made some observations on the minimum  degree of $(K_3, K_3)$-co-critical graphs and maximum number of possible edges of $(H_1, \ldots,  H_r)$-co-critical graphs.  Hanson and Toft~\cite{Hanson1987} in 1987 also studied  the minimum and maximum number  of edges over all $(H_1, \ldots,  H_k)$-co-critical graphs on $n$ vertices when $H_1, \ldots,  H_k$ are complete graphs, under the name of \dfn{strongly $(|H_1|, \ldots, |H_r|)$-saturated} graphs. Recently, this topic has  been studied under the name of \dfn{$\mathcal{R}_{\min}(H_1, \dots, H_k)$-saturated} graphs \cite{Chen2011, Ferrara2014,  RolekSong, SongZhang21}. We refer the reader to   a recent paper  by  Zhang and the second 
author~\cite{SongZhang21} for
further background on $(H_1, \ldots,  H_k)$-co-critical graphs.  Hanson and Toft~\cite{Hanson1987}    observed that for all $n\ge r(K_{t_1}, \dots, K_{t_k})$, the  graph $K_{r-2}+\overline{K}_{n-r+2}$ is $(K_{t_1}, \dots, K_{t_r})$-co-critical with   $ (r- 2)(n - r + 2) + \binom{r - 2}{2}$ edges.  They further  made the following conjecture that no   $(K_{t_1}, \dots, K_{t_r})$-co-critical graph on $n$ vertices can have fewer than $e(K_{r-2}+\overline{K}_{n-r+2})$ edges.

\begin{conj}[Hanson and Toft~\cite{Hanson1987}]\label{HTC}  Let   $r = r(K_{t_1}, \dots, K_{t_k})$. Then  
 every   $(K_{t_1}, \dots, K_{t_k})$-co-critical graph   on $n$ vertices has at least    
$$   (r- 2)(n - r + 2) + \binom{r - 2}{2}$$
edges. This bound is best possible for every $n$. 
\end{conj}
\medskip

It was shown in \cite{Chen2011} that every   $(K_3,  K_3)$-co-critical graph on $n\ge 56$ vertices has at least $4n-10$ edges, thereby verifying the first nontrivial case of Conjecture~\ref{HTC}. At this time, however, it seems that a complete resolution of Conjecture~\ref{HTC} remains elusive.  
 Some structural properties of $(K_3,  K_4)$-co-critical graphs are given in \cite{K3K4}.     Inspired  by Conjecture~\ref{HTC},  Ferrara,    Kim and  Yeager~\cite{Ferrara2014} proposed the following problem.
 
 \begin{prob}[Ferrara,    Kim and  Yeager~\cite{Ferrara2014}]\label{Problem} Let $H_1,   \ldots,  H_k$ be graphs, each with   at least one edge. Determine  the minimum number of edges of $(H_1,   \ldots,  H_k)$-co-critical graphs.
 \end{prob}
 
  In the same paper they settled Problem~\ref{Problem} when each $H_i$ is a matching of $m_i$ edges. 
  
  \begin{thm}[Ferrara,    Kim and  Yeager~\cite{Ferrara2014}]\label{matching}
Let $m_1,   \ldots, m_k$ be positive integers. Then every $(m_1K_2,   \ldots,  m_kK_2)$-co-critical graph on $n>3(m_1+\cdots m_k-k)$ vertices has at least $3(m_1+\cdots m_k)$ edges. This bound is best possible for all $n>3(m_1+\cdots m_k-k)$. 
\end{thm}
 
 Theorem~\ref{matching} yields the very first result  on Problem~\ref{Problem} for multicolor  $k$.   Very recently, Chen, Miao, Zhang and the second author~\cite{P3} settled  Problem~\ref{Problem}  when each $H_i=P_3$ by studying the minimum size of graphs $G$  such that  $\chi'(G)=\Delta(G)$ and    $\chi'(G+e)=\chi'(G)+1$ for  every  $e\in E(\overline{G})$, where $\chi'(H)$ denotes the chromatic index of a graph $H$.

 \begin{thm}[Chen, Miao, Song and Zhang~\cite{P3}]\label{P3k}
For all $k \ge 1$ and $n \ge k+1+ (k\text{ mod }2)$, every  $(P_3;k)$-co-critical graph $G$ on $n$ vertices satisfies  
 $$e(G)\ge {k \over 2}\left(n- \left\lceil {k \over 2} \right\rceil - \varepsilon\right) + {\lceil k/2 \rceil+\varepsilon \choose 2},$$
   where $\varepsilon=(n-\lceil k/2 \rceil)(\text{mod }  2)$. This bound is best possible for   all $k \ge 1$ and $n \ge \left\lceil {3k /2} \right\rceil +2$. 
\end{thm}

Motivated by Conjecture~\ref{HTC}, Rolek and the second author~\cite{RolekSong} recently initiated the study of the minimum number of possible edges over all $ (K_t, \mathcal{T}_k)$-co-critical graphs, where  $\mathcal{T}_k$ denotes  the family of all trees on $k$ vertices, and 
for all $t, k \ge 3$, we write $G\rightarrow (K_t, \mathcal{T}_k)$  if  for every $2$-coloring   $\tau: E(G) \to \{\text{red, blue} \}$, $G$ has either a red $K_t$ or a blue tree $T_k\in \mathcal{T}_k$; 
a non-complete graph $G$ is \dfn{$(K_t, \mathcal{T}_k)$-co-critical} if  
$G\nrightarrow (K_t, \mathcal{T}_k)$, but $G+e\rightarrow (K_t, \mathcal{T}_k)$ for all $e$ in $\overline{G}$.    
 Rolek and the second author~\cite{RolekSong} proved the following on $ (K_3, \mathcal{T}_k)$-co-critical graphs.

\begin{thm}[Rolek and Song~\cite{RolekSong}]\label{K3Tk}  Let $n, k\in \mathbb{N}$.  
\begin{enumerate}[(i)]
 \item  Every  $(K_3, \mathcal{T}_4)$-co-critical graph on $n\ge 18$ vertices has  at least $\left\lfloor 5n/2\right\rfloor$ edges. This bound is sharp for every $n\ge18$.  
 \item   For all  $k \ge 5$,  if $G$ is $(K_3, \mathcal{T}_k)$-co-critical on $n\ge  2k + (\lceil k/2 \rceil +1) \lceil k/2 \rceil -2$ vertices, then $$ e(G) \ge \left(\frac{3}{2}+\frac{1}{2}\left\lceil \frac{k}{2} \right\rceil\right)n-c(k),$$ 
 where $c(k)=\left(\frac{1}{2} \left\lceil \frac{k}{2} \right\rceil + \frac{3}{2} \right) k -2$. This bound is asymptotically best possible.
 \end{enumerate}
\end{thm}

 Very recently, Zhang and the second author~\cite{SongZhang21} obtained a lower bound for the size of $ (K_t, \mathcal{T}_k)$-co-critical graphs for all $t\ge 4$ and $k\ge\max\{6, t\}$. In addition, this bound is asymptotically best possible    when $t\in\{4,5\}$ and all $k\ge6$ and $n$ large.  They believe the lower bound is asymptotically best possible for all such $t$ and $k$.
 
\begin{thm}[Song and Zhang~\cite{SongZhang21}]\label{SongZhang}   
  Let  $ t, k\in \mathbb{N}$ with $t \ge 4$ and  $k\ge\max\{6, t\}$.  There exists a  constant  $\ell(t, k)$  such that, for  all $n\in \mathbb{N}$ with    $n \ge (t-1)(k-1)+1$,    if $G$  is  a $(K_t, \mathcal{T}_k)$-co-critical graph  on $n$ vertices, then     
  $$   e(G)\ge   \left(\frac{4t-9}{2}+\frac{1}{2}\left\lceil \frac{k}{2} \right\rceil\right)n-\ell(t, k).$$ 
This bound is asymptotically best possible    when $t\in\{4,5\}$ and all $k\ge6$ and $n \ge (2t-3)(k-1)+{\lceil k/2 \rceil}{\lceil k/2 \rceil}-1$.  
 \end{thm}

 The methods developed in \cite{RolekSong, SongZhang21}  may shed some light on attacking Conjecture~\ref{HTC}. In this paper, we continue to study the minimum number of possible edges over all  $(K_t, K_{1, k})$-co-critical graphs. 
By a classic result of Chv\'atal~\cite{Chvatal}, $r(K_t, K_{1,k})=(t-1)k+1$.  Hence,   every $(K_t, K_{1, k})$-co-critical graph has at least $r(K_t, K_{1,k})=(t-1)k+1$   vertices.     We prove the following main result. 
 
\begin{thm}\label{t:main}
For all $t\ge3$ and $k\ge3$, there exists a constant $\ell(t, k)$ such that  for all $n\in \mathbb{N}$ with  $n \ge (t-1)k+1$, if $G$ is a 
  $(K_t,K_{1,k})$-co-critical graph on  $n $ vertices, then    \[e(G)\geq \left(2t-4+\frac{k-1}{2}\right)n-\ell(t,k).\]
   This bound is asymptotically best possible    when $t\in\{3,4,5\}$ and all $k\ge3$ and $n \ge (2t-2)k+1$. 
 \end{thm}

We believe the bound in \cref{t:main} is  asymptotically best possible  for all $t\ge3$,    $k \ge 3$ and $n\ge (t-1)k+1$. It seems non-trivial to construct such extremal $(K_t, K_{1,k})$-co-critical graphs when $t\ge 6$.  Meanwhile, we  are able to obtain the   sharp bound  on the minimum number of edges of      $(K_3, K_{1,3})$-co-critical graphs on $n\ge 13$ vertices.

\begin{restatable}{thm}{threeclaw}\label{t:threeclaw}
 Every $(K_3, K_{1,3})$-co-critical graph on $n\ge13$ vertices  has at least   $  3n-4$ edges. This bound is sharp for all $n\ge13$.
\end{restatable}

To prove our main results, we first  establish some important structural properties of  $(K_t, K_{1,k})$-co-critical graphs in Section~\ref{property}. We then prove \cref{t:main} in Section~\ref{s:Ktstar} and \cref{t:threeclaw} in Section~\ref{s:K3claw}.

  \section{Structural properties of  $(K_t, K_{1,k})$-co-critical graphs}\label{property}
 
 We need to  introduce more notation. Given a graph $H$, a graph $G$ is \emph{$H$-free} if $G$ does not contain $H$ as a subgraph;  $G$ is \emph{$H$-saturated} if $G$ is $H$-free but $G+e$ is not $H$-free for every $e\in \overline{G}$.  
For a $ (K_t, K_{1,k})$-co-critical graph $G$,  let   $\tau : E(G) \to \{\text{red, blue} \}$  be a  $2$-coloring of $E(G)$ and let  $E_r$ and $E_b$ be the  color classes of the coloring $\tau$.  We use $G_{r}$ and $G_{b}$ to denote the spanning subgraphs of $G$  with edge sets  $E_r$ and $E_b$, respectively.  
We define  $\tau$ to be a  \emph{critical coloring} of $G$ if $G$ has neither  a red $K_t$ nor a blue $ K_{1,k}$ under $\tau$, that is,  if  $G_r$ is $K_t$-free and $G_b$ is $K_{1,k}$-free.  For every $v\in V(G)$, we use $d_r(v)$ and $N_r(v)$ to denote the degree and neighborhood of $v$ in $G_r$, respectively. Similarly, we define $d_b(v)$ and $N_b(v)$ to be the degree and neighborhood of $v$ in $G_b$, respectively.  \medskip   \medskip

Given  a $ (K_t, K_{1,k})$-co-critical graph $G$, it is worth noting  that   $G$ admits at least one critical coloring but,  for any edge $e\in E(\overline{G})$, 
$G+e$ admits  no critical coloring.  In order to find $e(G)$, we let $\tau : E(G) \to \{\text{red, blue} \}$  be a  critical coloring of $G$ such that $|E_r|$ is maximum among all critical colorings of $G$. Then $G_r$ is  $K_t$-free but, $G_r+e$ has a copy of $K_t$ for each $e\in E(\overline{G})$, that is, $G_r$ is $K_t$-saturated.
We state below two results on $K_t$-saturated graphs that shall be used in our proofs. We refer the reader to the  the dynamic survey \cite{Faudree2011} on  the extensive studies on $K_t$-saturated graphs.

\begin{thm}[Hajnal  \cite{Hajnal}]\label{t:Hajnal}
Let $t, n \in\mathbb{N}$. Let $G$ be a $K_t$-saturated graph on $n$ vertices.  Then either $\Delta(G) = n-1$ or $\delta(G) \ge 2(t-2)$. 
\end{thm}

\begin{thm}[Day \cite{Day2017}]\label{t:Day} 
Let $q \in\mathbb{N}$. There exists a constant $c = c(q)$ such that, for  all $3\le t \in \mathbb{N}$ and all $n\in \mathbb{N}$, if $G$ is a $K_t$-saturated graph on $n$ vertices with  $\delta(G) \ge q$, then $e(G) \ge qn - c$.
\end{thm}

 Following the ideas in \cite[Lemma 1.6]{RolekSong} and \cite[Theorem 7]{SongZhang21}, we   are now ready to prove a useful lemma on $(K_t, K_{1,k })$-co-critical graphs.

\begin{lem}\label{l:blue}
For all integers $t\ge3$ and  $k\ge3$, let $G$ be a  $(K_t, K_{1,k })$-co-critical graph on  $n\ge (t-1)k+1 $   vertices.  Among all critical colorings of $G$, let    $\tau : E(G) \rightarrow \{$red, blue$\}$ be  a critical coloring of $G$ with $|E_r|$  maximum. Let   $S\subseteq V(G)$ be  such that each vertex of $S$ has     degree at most $k-2$ in  $G_b$.   Then the following hold. 

\begin{enumerate}[(i)]

\item\label{alpha}   $ S $ is a clique in $G$, that is, $G[S]=K_{|S|}$.    Moreover, $\alpha(G_b[S])\le t-1$ and $|S|\le (t-1)(k-1)$.   

\item\label{n-2}   $\Delta(G_r) \le n-2$ and so   $\delta(G_r) \ge 2(t-2)$.
\item\label{triangle}  If $t=3$, then every  $e\in E_r$   belongs to at  most  $2k - 2$ triangles in $G$.
\item\label{n-3}  If $t=3$, then $G_r$ is $2$-connected and  $\Delta(G_r) \le n-3$. 
\end{enumerate}
\end{lem}

\begin{proof} 
To prove (\ref{alpha}), suppose   there exist $u,v\in S$ such that $uv\notin E(G)$.   Note that  $d_b(u)\leq k-2$ and $d_b(v)\leq k-2$ by the choice of $S$.   
   But then we obtain a  critical coloring of  $G+uv$ from $\tau$ by  coloring the edge $uv$ blue, contrary to the fact that  $G$ is $(K_t, K_{1,k })$-co-critical.  Hence $S$ is a clique in $G$. Since $G_r$ is $K_t$-free, it follows that $\alpha(G_b[S])\le t-1$ and so $|S|\le (t-1)(k-1)$.\medskip

To prove (\ref{n-2}),   
suppose   $\Delta(G_r) > n-2$. Let  $x \in V(G)$ with $d_{r}(x) =n-1$.   Note that $G_r\less x$ is $K_{t-1}$-free because $G_r$ is $K_t$-free. Since $G\ne K_n$,  there must exist  $y,z \in N_r(x)$ such that  $yz \not\in E(G)$.   But then we obtain a critical coloring  of $G+yz$ from  $\tau$ by first coloring the edge $yz$ red, and then recoloring $xy$ blue and all   edges incident with $y$ in $G_b$ red, contrary to the fact that $G+yz$ has no critical coloring. This proves that $\Delta(G_r) \le n-2$. Since $G_r$ is $K_t$-saturated, by \cref{t:Hajnal}, $\delta(G_r) \ge 2(t-2)$.  \medskip

To prove (\ref{triangle}), suppose there exists an edge $ uv \in E_r$ such that $uv$ belongs to at least $2k-1$ triangles in $G$.  
Since $G_r$ is  $K_3$-free,  we see that either $d_b(u)\ge k $ or $d_b(v)\ge k $. In either case,   $G_b$ contains a copy of $K_{1, {k}}$, a contradiction. \medskip

It remains to prove (\ref{n-3}).  
Suppose  first  $\Delta(G_r) \ge n-2$. By (iii), $\Delta(G_r) = n-2$. Let  $x,y \in V(G)$  be such that   $d_{r}(x) = n-2$ and  $y$ is the  unique non-neighbor of $x$ in $G_r$.  Note that   $N_r(x)$ is an independent set in $G_r$.  
    Suppose there exists a vertex  $z\in N_r(x)$  such that $yz$ is colored blue or $yz\notin E(G)$. In the former case, recoloring $yz$ red yields a  critical coloring of $G$ with $|E_r|+1$ red edges, contrary to the    choice of $\tau$. In the latter case, we obtain a critical coloring of $G+yz$ from $\tau$ by coloring the edge $yz$ red, a contradiction. This proves that   $y$ must be  complete to $N_r(x)$ in $G_r$. 
  Let $u\in N_r(x)$. Since $d_b(u)\le k-1$ and $|N_r(x)|=n-2\ge 2k-1$,  there exists a vertex $w\in N_r(x)$ such that $uw\notin E_b$.  Then $uw\notin E(G)$ because $N_r(x)$ is an independent set in $G_r$.  But then we obtain a critical coloring of  $G+uw$  from  $\tau$ by   first recoloring    edges  $ux, uy$ blue and  all  edges incident with $u$ in $G_b$ red, and then coloring the edge $uw$ red,    a contradiction.  Thus  $\Delta(G_r) \le n-3$. \medskip

Finally, we prove  that $G_r$ is $2$-connected. 
Suppose that $G_r$ is not $2$-connected. Since $G_{r}$ is $K_3$-saturated, we see that  $G_r$ is connected and  must have  a cut vertex, say $u$. 
Let $G_1$ and $G_2$ be two components of $G_r\setminus u$  such that $u$ is not complete to $V(G_2)$ in $G_r$.  
Let $v \in V(G_1)$ and $w\in V(G_1)\less N_r(u)$. Then
 $vw\in E(G)$, else we obtain a critical coloring of  $G+vw$ from $\tau$ by coloring the edge  $vw$ red. Thus  $vw\in E_b$. But then  we obtain a critical coloring of $G$ from $\tau$ by recoloring the   edge $vw$ red, contrary to the maximality of $|E_r|$. Therefore,  $G_r$ is $2$-connected. \medskip

This completes the proof of \cref{l:blue}. 
\end{proof}

\section{Proof of \cref{t:main}}\label{s:Ktstar}
\cref{t:main} follows immediately from \cref{t:Ktstar} and \cref{t:upper}. 

\begin{thm}\label{t:Ktstar}
For all $t\ge3$ and $k\ge3$, there exists a constant $\ell(t, k)$ such that  for all   $n\ge (t-1)k+1$, if $G$ is a 
  $(K_t,K_{1,k})$-co-critical graph on  $n $ vertices, then    \[e(G)\geq \left(2t-4+\frac{k-1}{2}\right)n-\ell(t,k).\]
  Moreover, $\ell(3, k)=(k-1)^2+5$ for all $k\ge3$.
 \end{thm}

\begin{proof}  Let $G$ be a  $(K_t, K_{1,k})$-co-critical graph on  $n \ge (t-1)k+1$   vertices, where  $t \ge 3$ and  $k \ge 3$.  Then      $G$ admits  a critical coloring. Among all critical colorings of $G$, let $\tau: E(G) \rightarrow \{$red, blue$\}$ be a critical coloring  of $G$ with $|E_r|$   maximum. By the choice of $\tau$, $G_r$ is $K_t$-saturated and $G_b$ is  $ K_{1,k}$-free.   Let   $S\subseteq V(G)$  be such that each vertex of $S$ has     degree at most $k-2$ in  $G_b$. By \cref{l:blue}(\ref{alpha}),   $|S|\le (t-1)(k-1)$, and so  \[e(G_b)\ge \frac{(k-1)(n-|S|)}2\ge \frac{(k-1)(n-(t-1)(k-1))}2=\frac{(k-1)n}2-\frac{(t-1)(k-1)^2}2.\]
We next estimate $e(G_r)$. By  \cref{l:blue}(\ref{n-2}),     $\delta(G_r)\ge 2t-4$.   Assume first $t\ge4$. By \cref{t:Day}, $e(G_r)\ge (2t-4)n-C(t)$ for some constant $C(t)$. 
Consequently,   
 \[e(G)=e(G_r)+e(G_b)\ge \left(2t-4+\frac{k-1}{2}\right)n-\ell(t,k),\]
where $\ell(t,k)=C(t)+(t-1)(k-1)^2/2$, as desired. \medskip

Assume next $t=3$. Then $\delta(G_r)\ge2$. We claim that $e(G_r)\ge 2n-5$. This is trivially true when $\delta(G)\ge4$. We may assume that $2\le \delta(G)\le3$. By \cref{t:p=3} and \cref{l:K3-sat}, $e(G_r)\ge 2n-5$, as claimed. 
Consequently,   
 \[e(G)=e(G_r)+e(G_b)\ge (2n-5)+\frac{(k-1)n}2- (k-1)^2= \left(2+\frac{k-1}{2}\right)n-\ell(3,k),\]
where $\ell(3,k)= (k-1)^2+5$, as desired. 
\end{proof}
 \bigskip
 
 As mentioned in the Introduction, we believe the bound in \cref{t:Ktstar} is  asymptotically best possible  for all $t\ge3$,   $k \ge 3$ and $n\ge(t-1)k+1$. However, it seems non-trivial to construct such extremal $(K_t, K_{1,k})$-co-critical graphs. We next show that \cref{t:Ktstar} is  asymptotically best possible  for $t \in \{3, 4, 5\}$,   all $k \ge 3$ and $n\ge (2t-2)k+1$.   \medskip

\begin{thm}\label{t:upper}
For each $t \in \{3, 4, 5\}$,  all $k \ge 3$ and  $n \ge (2t-2)k+1 $, there exists a   $(K_t, K_{1,k})$-co-critical graph  $G$ on $n$ vertices such that 
$$e(G) = \left(2t-4+\frac{k-1}{2}\right)n+C(t, k),$$
where $C(3, k)=k^2 - 3k - 3-\frac{(k-1)\cdot \eps}2$,  $C(t, k)=\frac{1}{2}(k-1)\big((t^2-t-2)k-(t^2+t+\eps-6)-(t-2)\big(k+\frac{1}{2}(5t-7)\big)$ for   each $t\in\{4,5\}$,  and    $\eps=1$ if  $k$ or $n$ is odd;   $0$ if both $k$ and $n$ are even. \end{thm}

\begin{proof}
Let $t, k, n,\eps  $ be as given in the statement. 
By the choice of $\eps$, we see that $n-(2t-3)k-\eps\ge k $ and $(k-1)(n- (2t-3)k-\eps)$ is always even.   We shall construct a  $(K_t, K_{1,k})$-co-critical graph on $n\ge (2t-2)k+1$ vertices which yields the desired number of edges. \bigskip

We begin with the case $t=3$. Let $A: =K_{k-1}$,   $H=(B,C)$ be a $(k-2)$-regular bipartite graph with $|B|=|C|=k-1$,  and let   $R $ be a $(k-1)$-regular graph on $n-3k$ vertices when $\eps=0$, and   the disjoint union of $K_1$ and a  $(k-1)$-regular graph on $n - 3k -1 $ vertices when $\varepsilon=1$.       Let $G_1$ be obtained from disjoint copies of $A,  H, R$ by joining every vertex in $A$ to all vertices in $B$.  Finally, let $G$ be the graph obtained from $G_1$ by adding three  new vertices $x, y, z$, and then  joining: $x$ to  every vertex in $V(G_1)$;  $y$ to   every vertex in every vertex in $V(G_1)\cup\{z\}$;   $z$ to  every vertex in $V(A)\cup C$.     The construction of  $G$  is  depicted in Figure~\ref{uniqueGG}.    It can be easily checked that  $e(G) = \left(2+\frac{k-1}{2}\right)n + k^2 - 3k-3 -\frac{(k-1)\cdot \varepsilon}2$. It suffices to show that $G$ is  $(K_3, K_{1,k})$-co-critical.   \vspace{-0.5cm}

\begin{figure}[htbp]
 \centering
 \includegraphics[scale=1.8]{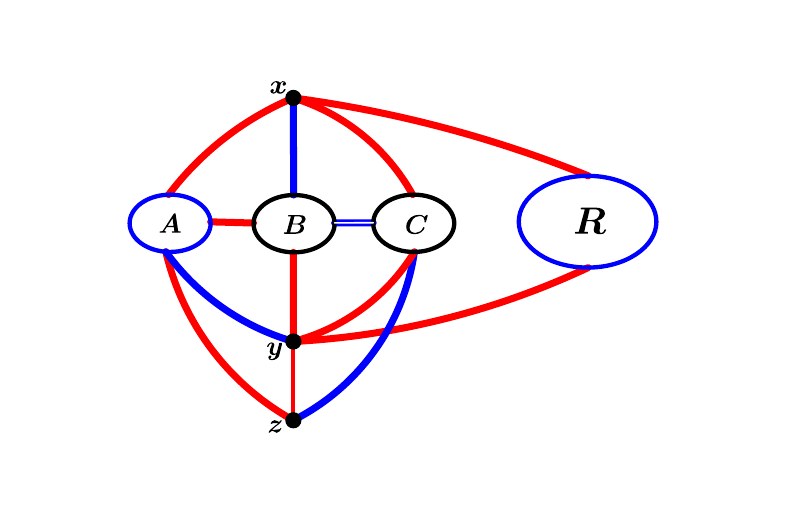}
 \vspace{-0.5cm}
   \caption{A $(K_3,K_{1,k})$-co-critical graph for all $k\ge3$  with a unique critical coloring.}
\label{uniqueGG} 
 \end{figure}

Let  $\sigma: E(G) \rightarrow \{$red, blue$\}$ be defined as follows: all edges in   $A, H$ and $R$ are colored blue;  all edges between $x $ and $B$ are colored blue, all edges between $y$ and $V(A)$ are colored blue, and  all edges between $z$ and $C$ are colored blue; the remaining edges of $G$ are all  colored red. The coloring    $\sigma $ is  depicted in Figure \ref{uniqueGG}.   It is simple to check that   $\sigma $   is a critical coloring  of $G$. We next show that $\sigma$ is the unique critical coloring  of $G$. \medskip

Let $\tau: E(G) \rightarrow \{$red, blue$\}$ be an arbitrary critical coloring of $G$.    Let $G^\tau_r$ and $G^\tau_b$ be  $G_r$ and $G_b$ under the coloring $\tau$, respectively. Then $G^\tau_r$ is $K_3$-free and $G^\tau_b$ is $K_{1,k}$-free. Thus $\Delta(G^\tau_b)\le k-1$.  Note that every edge of $A$  belongs to $2k-1$ triangles in $G$,   and must be colored blue under $\tau$ by \cref{l:blue}(\ref{triangle}). Next, since $\Delta(G^\tau_b)\le k-1$, we may assume that $yw\in E(G^\tau_r)$   for some $w\in B\cup\{z\}$.  We claim that $x$ is  complete to $V(A)$   in $G^\tau_r$.  Suppose   $xu\in E(G^\tau_b)$   for some $u\in V((A)$. Then $u$ is  complete to $B\cup \{y, z\}$ in $G^\tau_r$. But then  $G^\tau_r[\{y,u,w\}]=K_3$, a contradiction.  Thus, $x$ is complete to $V(A)$  in $G^\tau_r$, as claimed.  We next claim that  $y$ is complete to $V(A)$   in $G^\tau_b$. Suppose  $yv \in E(G^\tau_r)$ for some $v\in V(A)$. Then   $vw$ must be colored blue under $\tau$,   else $G^\tau_r\{[y,v,w]\}=K_3$. Hence $v$ is  complete to $(B\cup\{z\})\setminus\{w\}$ in $G^\tau_r$, and so $y$ is complete to $(B\cup\{z\})\setminus\{w\}$ in $G^\tau_b$. Furthermore, $y$ is complete to $C\cup V(A) $ in $G^\tau_r$. Then  $w$ is complete to $V(A)$ in $G^\tau_b$. Let $w'\in C$ be  such that $ww'\in E(G)$. Since $G_r[\{y, w, w'\}]\ne K_3$, we see that  $ww'\in  E(G^\tau_b)$. But then    $w$ is complete to $V(A)\cup\{w'\}$ in $G^\tau_b$, contrary to the fact that $G^\tau_b$ is $K_{1,k}$-free.  This proves that  $y$ is complete to $V(A)$   in $G^\tau_b$, as claimed.  It then follows that in $G^\tau_r$: $V(A)$ is  complete to $B\cup\{x, z\}$;   $y$ is  complete to $\{z\}\cup B\cup C\cup V(R)$. Thus, all edges  in $H\cup R$,  all edges between $x$ and $B$ and all edges between $z$ and $C$   must be colored blue under $\tau$.  This proves that  $ \tau=\sigma$. Therefore,   $\sigma$  is indeed the unique critical coloring of $G$. It is easy to check that for each pair $u,v\in V(G)$ with $uv\notin E(G)$, $G+uv$ contains a red $K_3$ if $uv$ is colored red and a blue $K_{1,k}$ if $uv$ is colored blue.  Hence $G$ is  a desired  $(K_3, K_{1,k})$-co-critical graph. \vspace{-0.5cm}

\begin{figure}[htb]
\centering
\includegraphics[scale=1.12]{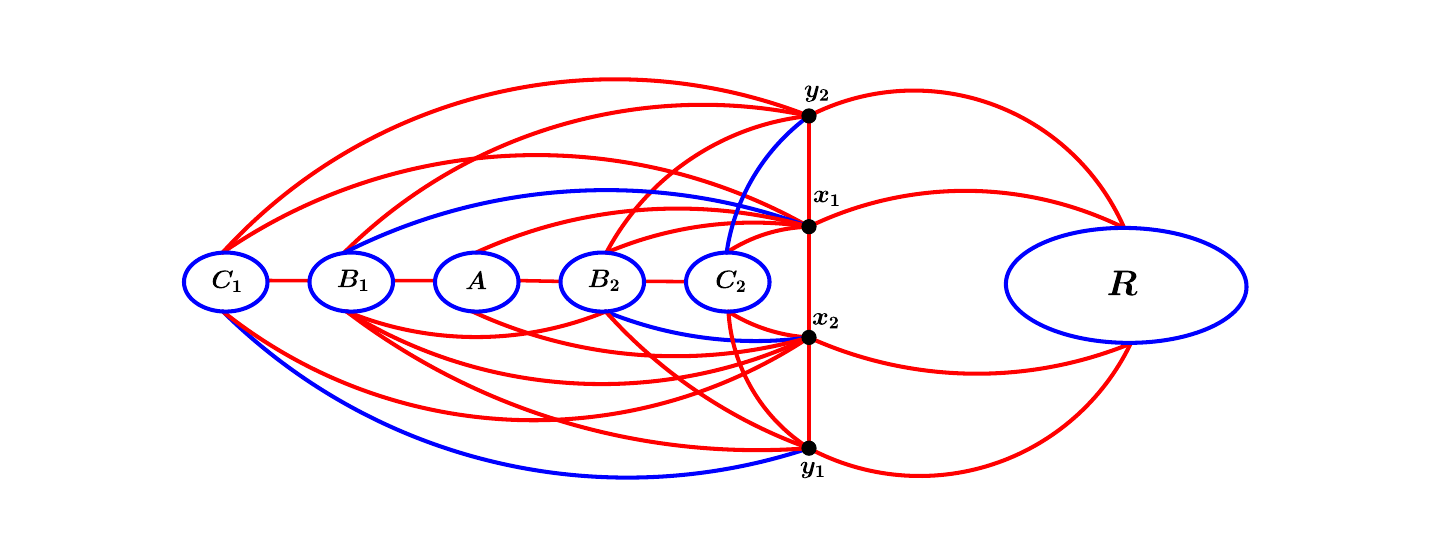}
\vspace{-0.5cm}
\caption{A  $(K_4, K_{1,k})$-co-critical  graph for all $k \ge 3$.}
\label{K4Tk-Saturated}
\end{figure}
\begin{figure}[htb]
\centering
\includegraphics[scale=1.02]{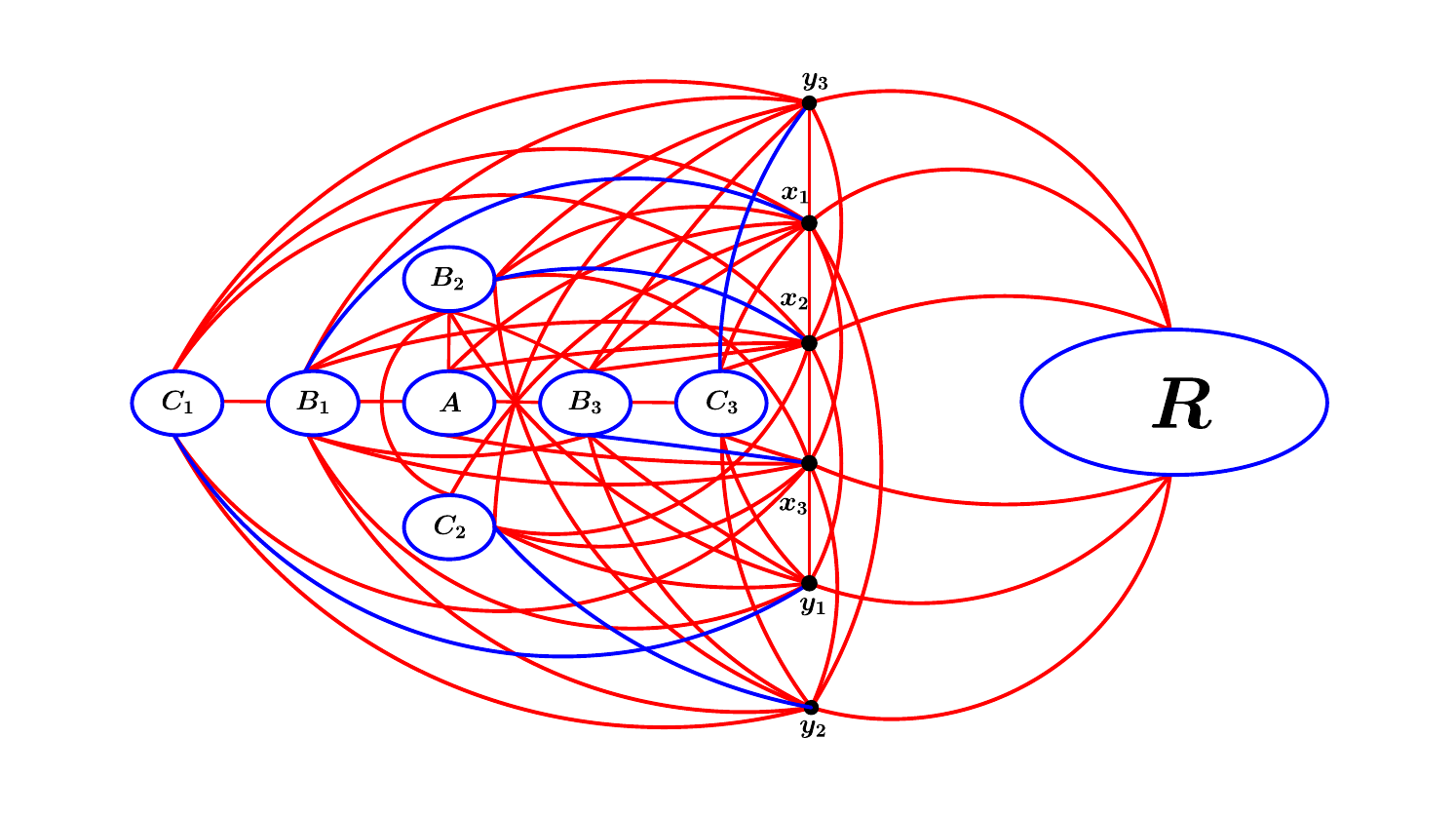}
\vspace{-0.5cm}
\caption{A  $(K_5, K_{1,k})$-co-critical graph for all $k \ge 3$.}
\label{K5Tk-Saturated}
\end{figure}

It remains to  consider the case $t\in\{4,5\}$.  Following the ideas in \cite[Theorem 9]{SongZhang21}, 
let $A: =K_{k}$. For each $i \in [t-2]$,  let  $B_i: =K_{k-1}$ and $C_i:=K_{k-1}$.   Let $H$ be obtained from disjoint copies of $A, B_1, \ldots, B_{t-2}, C_1, \ldots, C_{t-2}$ by joining every vertex in $B_i$ to all vertices in $A\cup C_i\cup  B_j$ for each   $i\in [t-2]$ and  all $  j \in [t-2]$ with $j \neq i$.  Let $R $ be a $(k-1)$-regular graph on $n-(2t-3)k$ vertices when $\eps=0$, and     the disjoint union of $K_1$ and a  $(k-1)$-regular graph on $n-(2t-3)k-1$ vertices when $\eps=1$.     Finally, let $G$ be the graph obtained from $H  \cup R$ by adding $2t-4$ new vertices $x_1, \ldots, x_{t-2}, y_1, \ldots, y_{t-2}$, and then, for each   $i\in [t-2]$, joining: $x_i$ to  every vertex in $V(H)\cup V(R)$ and all $x_j$;  and $y_i$ to  every vertex in $(V(H) \setminus V(A))\cup V(R)$  and all  $x_j$, where  $  j \in [t-2]$ with $j \neq i$. The construction of  $G$ when $t=4$ and $k \ge 3$ is depicted in Figure~\ref{K4Tk-Saturated}, and the construction of  $G$ when $t=5$ and $k \ge 3$ is depicted in Figure~\ref{K5Tk-Saturated}.  Let $X:=\{x_1,\ldots,x_{t-2}\}$ and $Y:= \{y_1, \ldots, y_{t-2}\}$.  Note that 
\begin{align*}
e_G(X\cup Y, V(G)\less (X\cup Y))&=(t-2)(n-(2t-4))+(t-2)(n-(2t-4+k))\\
&=(t-2)(2n-4t-k+8);
\end{align*}
 $e_G(X\cup Y)= {t-2 \choose 2}+(t-2)(t-3)$; $e_G(V(B_1) \cup \cdots \cup V(B_{t-2}), V(C_1) \cup \cdots \cup V(C_{t-2}))=(t-2)(k-1)^2$; $e_G(V(C_1) \cup \cdots \cup V(C_{t-2}))=(t-2) {k-1\choose 2}$; $e_G(V(A)\cup V(B_1) \cup \cdots \cup V(B_{t-2}))={(t-2)(k-1)+k  \choose 2}$; $e(H)=(k-1)(n-(2t-3)k-\eps)/2$.    Therefore,  
\begin{align*}
e(G)&= (t-2)(2n-4t-k+8)+{t-2 \choose 2}+(t-2)(t-3)+(t-2)(k-1)^2 +\\
&\quad (t-2) {k-1\choose 2}+{(t-2)(k-1)+k \choose 2} + (k-1)(n-(2t-3)k-\eps)/2\\
&=\left(2t-4+\frac{k-1}2\right)n-(t-2)k-\frac{1}{2}(t-2)(5t-7)\\
&\quad+(k-1) \big((t-2)(k-1)+(t-2)(k-2)/2+(t-1)(tk-k-t+2)/2 \big) \\
&\quad  -(k-1)((2t-3)k+\eps)/2\\
&=\left(2t-4+\frac{k-1}2\right)n-(t-2)\big(k+\frac{1}{2}(5t-7)\big)+\frac{1}{2}(k-1)\big((t^2-t-2)k-(t^2+t+\eps-6)\big) \\
&=\left(2t-4+\frac{k-1}2\right)n+C(t,k),
\end{align*}
where $C(t,k)=\frac{1}{2}(k-1)\big((t^2-t-2)k-(t^2+t+\eps-6)-(t-2)\big(k+\frac{1}{2}(5t-7)\big)$.  \medskip

It suffices to show that $G$ is $(K_t, K_{1,k})$-co-critical.  Let  $\sigma: E(G) \rightarrow \{$red, blue$\}$ be defined as follows: all edges in   $A, B_1, \ldots, B_{t-2}$, $C_1, \ldots, C_{t-2}$ and $R$ are colored blue; for every $i\in[t-2]$,  all edges between $x_i$ and $B_i$ are colored blue and  all edges between $y_i$ and $C_i$ are colored blue; the remaining edges of $G$ are all  colored red. The coloring    $\sigma $ is  depicted in   Figure \ref{K4Tk-Saturated} when $t=4$ and $k \ge 3$;  in   Figure \ref{K5Tk-Saturated} when $t=5$  and $k \ge 3$. It is simple to check that   $\sigma $   is a critical coloring  of $G$. We next show that $\sigma$ is the unique critical coloring  of $G$  up to symmetry.  \medskip

  Let $\tau: E(G) \rightarrow \{$red, blue$\}$ be an arbitrary critical coloring of $G$ with color classes $E_r$ and $E_b$.  It suffices to show that $\tau=\sigma$ up   to symmetry.  Let $G^\tau_r$ and $G^\tau_b$ be  $G_r$ and $G_b$ under the coloring $\tau$, respectively.  Let $W:=V(A) \cup V(B_1) \cup \cdots \cup V(B_{t-2}) \cup X$. Note that $G[W] = K_{(t-1)k}$ and $\alpha(G^\tau_b[W])\le t-1$.  We claim that $G^\tau_b[W]=(t-1)K_k$. Let $S$ be a maximum independent set of $G^\tau_b[W]$. Then $|S|\le t-1$ and   every vertex of  $  W\less S$  is adjacent to at least one vertex of  $S$ in $G^\tau_b[W]$. Since $G^\tau_b$ is $K_{1,k}$-free, we see that 
\[(t-1)(k-1)\le |W\less S|\le e_{G^\tau_b}(W\less S, S)= \sum_{v\in S} d_{G^\tau_b}(v)\le |S|(k-1)\le (t-1)(k-1).\]
It follows that $|S|=t-1$,   $d_{G^\tau_b}(v)=k-1$ for each $v\in S$, and every vertex of  $  W\less S$  is adjacent to exactly one vertex in $S$ in $G^\tau_b[W]$.  Therefore,  each vertex of $S$ belongs to a blue $K_k$  under $\tau$ because $\alpha(G^\tau_b[W])\le t-1$. Thus  $G^\tau_b[W]=(t-1)K_k$, as claimed.  It then  follows  that: for each $i  \in [t-2]$,  $y_i$ is complete to $ V(B_1) \cup \cdots \cup V(B_{t-2}) \cup (X\less x_i)$  in $G^\tau_r$;  and  $V(C_i)$ is complete to $V(B_i)\cup X$ in $G^\tau_r$.  Let $D_1, \ldots, D_{t-1}$ be the $t-1$ components of $G^\tau_b[W]$. We next prove several   claims. \bigskip

\setcounter{counter}{0}
\noindent {\bf Claim\refstepcounter{counter}\label{A}  \arabic{counter}.}  
  $ A= D_i $ for some $i\in [t-1]$. 

\begin{proof} Suppose $ A  \neq  D_i $ for all $i\in [t-1]$. Then  for each $i\in [t-1]$,  we see that  \[(V(B_1) \cup \cdots \cup V(B_{t-2})\cup X) \cap V(D_i) \neq \emptyset.\] Let $d_i \in (V(B_1) \cup \cdots \cup V(B_{t-2})\cup X) \cap V(D_i)$ for each $i\in [t-1]$. Then $d_1, \ldots, d_{t-1}$ are pairwise distinct and $G^\tau_r[\{d_1, \ldots, d_{t-1} \}]=K_{t-1}$.  Note that either $X\subseteq \{d_1, \ldots, d_{t-1}\}$ or $x_i\notin  \{d_1, \ldots, d_{t-1}\}$ for some $i\in[t-2]$.  It follows that:  in the former case, we may assume that $V(B_1)\cap  \{d_1, \ldots, d_{t-1}\}\ne \emptyset$, and then $G^\tau_r[\{d_1, \ldots, d_{t-1}, u\}]=K_t$ for any $u\in V(C_1)$;  in the latter case,   $G^\tau_r[\{d_1, \ldots, d_{t-1}, y_i\}]=K_t$, because  $y_i$ is complete to $V(B_1) \cup \cdots \cup V(B_{t-2})\cup (X\less x_i)$ in $G^\tau_r$. In both cases we obtain  a contradiction because $G_r^\tau$ is $K_t$-free.  \end{proof}\medskip

 By Claim~\ref{A}, we may assume that $ A = D_{t-1} $. Then $V(A)$ is complete to $ V(B_1) \cup \cdots \cup V(B_{t-2}) \cup X$ in $G^\tau_r$. Since $G^\tau_b$ is  $K_{1,k}$-free, we see that   \medskip
 
  \noindent {\bf Claim\refstepcounter{counter}\label{ci}  \arabic{counter}.} 
      for each $i \in [t-2]$, either $E_r\cap E(C_i)\ne\es$, or every vertex of $C_i$ is adjacent to   at least $t-3$ vertices  of   $Y$ in $G^\tau_r$. 
\medskip

\noindent {\bf Claim\refstepcounter{counter}\label{BD1}  \arabic{counter}.} 
   $|X \cap V(D_i)|= 1$ for each  $i \in [t-2]$. 
   
  \begin{proof} Suppose $|X \cap V(D_i)|\ne 1$ for some $i\in[t-2]$.  Since $|X|=t-2$, we may assume that  $|X \cap V(D_1)|\ge2$ and $X \cap V(D_{t-2})=\es$.  We may further assume that    $x_1, x_2\in V(D_1)$.  Then $x_1x_2\in E_b$.   Since $X \cap V(D_{t-2})=\es$ and  for all $i\in [t-2]$, $|V(B_i)|=k-1<k=|V(D_{t-2})|$, we may assume that $V(B_i)\cap V(D_{t-2})  \ne \es$  for each $i\in[2]$.  Let $b_i\in V(B_i)\cap V(D_{t-2})$ for each $i\in[2]$.   By Claim~\ref{ci},   for each $i\in[2]$, let   $c_i c_i^* \in E_r$, where $c_i\in V(C_i)$ and $c^*_i\in V(C_i)\cup Y$. Note that  $G^\tau_r[\{ c_1,  b_1, x_1, c_1^*\}]=K_4$ if $c_1^*\ne y_1$; and $G^\tau_r[\{ c_1,  b_1, x_2, c_1^*\}]=K_4$ if $c_1^*= y_1$. Thus  $t=5$. We next claim that 
  $V(B_i)\cap V(D_2)=\es$  for each $i\in[2]$.  Suppose, say $V(B_1)\cap V(D_2)\ne\es$. Let $d_1\in V(B_1)\cap V(D_2)$. Then  $G^\tau_r[ \{c_1, b_1, d_1, x_1, c_1^*  \}]=K_5$ if  $c_1^*\ne y_1$; and $G^\tau_r[ \{c_1, b_1, d_1,   x_2, c_1^*  \}]=K_5$ if  $c_1^*= y_1$,  a contradiction. Thus $V(B_i)\cap V(D_2)=\es$  for each $i\in[2]$, as claimed. Then $D_2=B_3\cup\{x_3\}$. Note that    $G^\tau_r[\{c_1, b_1,   x_1, x_3, c_1^*\}]=K_5$  if  $c_1^*\in V(C_1)$. It follows that  $c_1^*\in Y$  and $E_r\cap E(C_1)=\es$.  By Claim~\ref{ci},  every vertex of $C_1$ is adjacent to at least two vertices of $Y$ in $G_r^\tau$. We may further assume that $c_1^*= y_i$ for some $i\in[2]$. But then   $G^\tau_r[ \{c_1, b_1,   x_{3-i},  x_3, y_i\}]=K_5$, a contradiction.   
 \end{proof}\medskip    
 
\noindent {\bf Claim\refstepcounter{counter}\label{BD}  \arabic{counter}.} 
   For each  $i \in [t-2]$, $V(B_i) \subseteq V(D_j)$ for some  $j \in [t-2]$. 

\begin{proof} Suppose there exists an $i\in [t-2]$ such that $V(B_i) \nsubseteq V(D_j)$ for every  $j \in [t-2]$. We may assume $i=1$. By Claim~\ref{ci},   let   $c_1 c_1^* \in E_r$, where $c_1\in V(C_1)$ and $c^*_1\in V(C_1)\cup Y$.   
We claim that $V(B_1)\cap V(D_j) =\es$ for some $j\in [t-2]$. Suppose 
$V(B_1)\cap V(D_j)\ne\es$ for all $j\in [t-2]$.
 Let   $d_j \in V(B_1) \cap V(D_j)$ for all $j\in [t-2]$.  But then $G^\tau_r[\{d_1,  \ldots, d_{t-2}, c_1, c_1^*\}] =K_t$, a contradiction.  Thus $V(B_1)\cap V(D_j) =\es$ for some $j\in [t-2]$, as claimed. We may assume that $V(B_1)\cap V(D_{t-2}) =\es$. Since $V(B_1) \nsubseteq V(D_j)$ for every  $j \in [t-2]$, it follows that  $t=5$,    $V(B_1) \subseteq V(D_1) \cup   V(D_{2})$,  and  $V(B_1)\cap V(D_{1})\ne\es$  and $V(B_1)\cap V(D_{2})\ne\es$.   Let $d_1 \in V(B_1) \cap V(D_1)$ and $d_2 \in V(B_1) \cap V(D_2)$.  By Claim~\ref{BD1}, let $x_i\in X\cap V(D_3)$.   Then  $G^\tau_r[\{d_1, d_2, x_i, c_1, c_1^*\}]=K_5$ if  $c^*_1\in V(C_1)$. Thus $c^*_1\in    Y$ and    $E_r\cap E(C_1)=\es$. By Claim~\ref{ci}, $c_1$ is adjacent to at least two vertices in $Y$ in $G_r^\tau$. We may assume that $c^*_1\ne y_i$. Then $G^\tau_r[\{d_1, d_2, x_i, c_1, c_1^*\}]=K_5$, 
  a contradiction.   
   \end{proof}\medskip

By  Claim \ref{BD1} and  Claim \ref{BD},    $V(B_i) \cup V(B_j) \nsubseteq D_{\ell}$ for any $i \ne j \in [t-2]$ and all $\ell \in [t-2]$. 
  By symmetry,    we may assume that $V(B_i) \subseteq V(D_i)$ for all $i \in [t-2]$. Then $V(B_i) \cup \{x_j\} = V(D_{i})$ for some $j  \in [t-2]$ since $|V(D_{i})|=|V(B_i)|+1$ and $V(B_1) \cup \cdots \cup V(B_{t-2})\cup X=V(D_1) \cup \cdots \cup V(D_{t-2})$. 
By symmetry, we may assume that $V(B_i) \cup \{x_i\} = V(D_{i})$ for all $i \in [t-2]$. It follows that for all $i, j \in [t-2]$ with $i \neq j$, $V(B_i)$ is complete to $V(B_j)$ in $G^\tau_r$, $x_i$ is complete to $X \less x_i$ and $V(B_j)$ in $G^\tau_r$, $y_i$ is complete to $C_i$ in $G^\tau_b$, $y_i$ is complete to $V(C_j) \cup (X \less x_i)$ in $G^\tau_r$, $x_i$ is complete to $V(B_i)$ in $G^\tau_b$, $\{x_i, y_i\}$ is complete to $V(R)$ in $G^\tau_r$, all edges in $A, B_1, \ldots, B_{t-2}, C_1, \ldots, C_{t-2}$ and $R$ are colored blue under $\tau$. This proves that $\tau=\sigma$ and thus $\sigma$ is the unique critical coloring of $G$ up  to symmetry. \medskip

It can be easily checked that adding any edge $e\in E(\overline{G})$ to $G$ creates a red $K_t$ if $e$ is colored red, and a blue $K_{1,k}$ if $e$ is colored blue. Hence, $G$ is $(K_t, K_{1,k})$-co-critical, as desired.  \medskip

 This completes the proof of  \cref{t:upper}. \end{proof}

\section{Proof of \cref{t:threeclaw}}\label{s:K3claw}

In this section, we prove \cref{t:threeclaw}.    We shall need   a   result  of Duffus and Hanson~\cite{Duffus1986} on $K_3$-saturated graphs with minimum degree   $ 3$, and  a   result   of Rolek and the second author~\cite{RolekSong} on $K_3$-saturated graphs with minimum degree   $ 2$. 

\begin{thm}[Duffus and Hanson~\cite{Duffus1986}]\label{t:p=3}

Let  $G$ be a $K_3$-saturated graph on  $n  \ge 10$ vertices with     $\delta(G) = 3$.  Then $e(G) \ge  3n - 15$. 
\end{thm}

\begin{figure}[htbp]
\centering
 \includegraphics[scale=0.5]{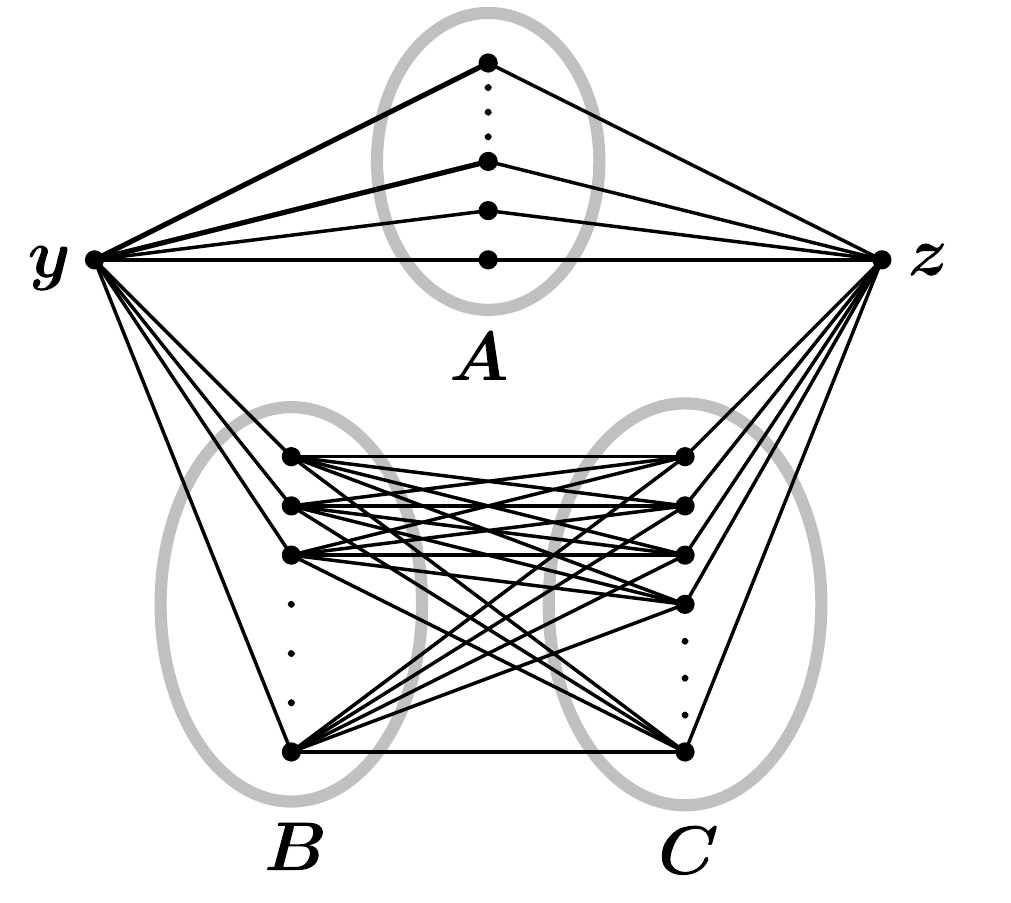} 
   \caption{Graph $J$.}
  \label{J}
  \end{figure}

The graph $J$ depicted  in  Figure~\ref{J}  is a $K_3$-saturated graph with $|J|\ge5$ and   $\delta(J)=2$, where $A\ne\emptyset$ and  either  $B=C=\emptyset$ or  $B\ne\emptyset$ and  $C \ne \emptyset$;   $A$, $B$ and  $C$ are  independent sets in $J$ and pairwise disjoint;    $A$ is anti-complete to $B\cup C$ and  $B$ is complete to $C$;  $N_J(y)=A\cup B$ and $N_J(z)=A\cup C$;  and $|A|+|B|+|C|=|J|-2$.  It is straightforward to check that 
$e(J) = 2(|J| - 2) + |B||C| - |B| - |C|\ge 2|J|-5$. Moreover, $e(J)=2|J|-5$ when $|B|=1$ or $|C|=1$. That is, $e(J)=2|J|-5$ when  $J$ is  obtained from $C_5$ by   repeatedly duplicating vertices of degree $2$.

\begin{lem}[Rolek and Song~\cite{RolekSong}] \label{l:K3-sat}
Let $G$ be a $K_3$-saturated graph on $n\geq 5$ vertices with $\delta(G) =2$. Then
   $G \cong J$ and $e(J)\ge 2n-5$. Furthermore,   if  $e(G) = 2n-k$ for some $k\in\{0,1,2,3,4,5\}$,   then   $|B||C| - |B| - |C| = 4-k$, where $A,B,C,$ and $J$ are as depicted in Figure~\ref{J}  and the values of $|B|$ and $|C|$ are summarized as below:

\begin{center}
\begin{tabular}{l l l}
\hline
$k$ & $e(J)$ & values of $|B|$ and $|C|$ with $|B|\leq |C|$ \\\hline
5 & $2n-5$ & $ |B|=1$ and $|C|\geq 1$  \\
4 & $2n-4$ & $|B|=|C|=2$ or $|B| = |C| = 0$ \\
3 & $2n-3$ & $|B|=2$ and $|C| = 3$ \\
2 & $2n-2$ & $|B|=2$ and $|C| = 4$ \\
1 & $2n-1$ & $|B|=2$ and $|C|= 5$ or $|B|=|C|=3$ \\
0 & $2n$ & $|B|=2$ and $|C|=6$ \\\hline
\end{tabular}
\end{center}
\end{lem}
\bigskip

We are now ready to prove  \cref{t:threeclaw}, which we restate here for convenience.\threeclaw*

\begin{proof} The sharpness of the bound for each $n\ge13$ is witnessed by \cref{t:upper} when $t=k=3$ (and so $\eps=1$).  
We next show that every   $(K_3,K_{1,3})$-co-critical graph  on $n\geq 13$ vertices has at least  $  3n-4$ edges.    Suppose this is not true. Let  $G$ be a  $(K_3,K_{1,3})$-co-critical graph on $n\geq 13$ vertices with $e(G)\leq 3n-5$. Then
   $G$ has at least one critical coloring. Among all critical colorings of $G$, let $\tau:E(G)\longrightarrow \{$red, blue$\}$ be a critical coloring of $G$ with $|E_r|$ maximum.  By the choice of $\tau$, $G_r$ is $K_3$-saturated and $G_b$ is $K_{1,3}$-free.  Thus $\Delta(G_b)\leq 2$. 
 We next prove several claims. \\

\setcounter{counter}{0} 

\noindent {\bf Claim\refstepcounter{counter}\label{Gb}  \arabic{counter}.}   
$n-2 \leq e(G_b) \leq n$   and  $e(G_r)\leq 2n-3$. 

\begin{proof}      Since  $\Delta(G_b)\leq 2$, we see that   $e(G_b)\leq n$. 
Suppose $e(G_b)\le n-3$. Then $\delta(G_b)\le 1$.  Let   $S\subseteq V(G)$ be  such that each vertex of $S$ has     degree at most one in  $G_b$. By \cref{l:blue}(\ref{alpha}), $S$ is a clique in $G$,  $\alpha(G_b[S])\le2$ and  $|S|\le 4$. Thus   $G_b[S]$ has at most two vertices of degree zero.   Note that 
 \[\sum_{v\in S} d_b(v)+2(n-|S|)=2e(G_b)\le 2n-6.\] Thus $|S|\ge 
3+\sum_{v\in S} d_b(v)/2$. It follows that $|S|=4$ and $e(G_b[S])=1$, which is impossible because $\alpha(G_b[S])\le2$. This proves that   $e(G_b)\ge n-2$ and so    $e(G_r)=e(G)-e(G_b)\leq 2n-3$.   
\end{proof}    

 \noindent {\bf Claim\refstepcounter{counter}\label{Gr}  \arabic{counter}.}
$\delta(G_r)=2$, and so $G_r \cong J$, where    $A, B, C, y,z$ are as pictured in Figure \ref{J}.  
 
\begin{proof} By  \cref{l:blue}(\ref{n-3}),   $G_r$ is 2-connected and so   $\delta(G_r)\geq 2$. Suppose $\delta(G_r) \ge3$.
By Claim~\ref{Gb},  $e(G_r)\le 2n-3$ and so  $\delta(G_r)\le3$. Thus  $\delta(G_r)=3$.  By \cref{t:p=3} and  Claim~\ref{Gb},    $e(G_r)\ge 3n-15$  and $e(G_b)\ge n-2$. But then $e(G) = e(G_r) + e(G_b) \geq 4n-17 > 3n-5$ because  $n\geq 13$, a contradiction. This proves that  $\delta(G_r)=2$. By  \cref{l:K3-sat},  $G_r \cong J$,  where    $A, B, C, y,z$ are as pictured in Figure \ref{J}.  
\end{proof} 

 By Claim~\ref{Gr}, $G_r\cong J$,  where    $A, B, C, y,z$ are as pictured in Figure \ref{J}.  
  By the definition of the graph $J$, 
we see that   $A\ne\emptyset$  and  either  $B=C=\emptyset$ or  $B\ne\emptyset$ and  $C \ne \emptyset$; $A, B, C$ are pairwise independent sets in $J\less\{y,z\}$;       $yz\notin E(J)$; and $N_J(y)=A\cup B$ and $N_J(z)=A\cup C$.    We may assume that  $|B|\leq |C|$. If $|B|=0$, then  $|B|=|C|=0$. But then  $\Delta(G_r)=d_r(y) = d_r(z) = n-2$, contrary to  the fact that $\Delta(G_r)\le n-3$ by \cref{l:blue}(\ref{n-3}). Thus $|C|\ge |B|\ge1$. Moreover,     \cref{l:K3-sat} implies that  $|B|\le2$ because   $e(G_r)\le 2n-3$. Let $B := \{z_1, \ldots, z_{|B|}\}$ and $C := \{y_1, \ldots, y_{|C|}\}$. \\ 
 
\noindent  {\bf Claim\refstepcounter{counter}\label{B}  \arabic{counter}.}
 $B=\{z_1, z_2\}$, $2\le |C|\le3$ and  $e(G_r) \ge2n-4$.

\begin{proof}   
Suppose   $B=\{z_1\}$.  By symmetry of $A$ and $C$, we may assume that $|A|\ge|C|$. Then  $|A|\ge5$ because $n\ge13$. Since $d_b(z_1)\le 2$, there must exist a vertex  $v\in A$   such that $z_1 v\notin E(G)$. But then we  obtain a critical coloring of $G+z_1v$ from $\tau$ by:  first  recoloring all edges incident with $y$ or $z_1$ in $G_b$ red,   then  recoloring  $yz_1  $ blue, and   coloring  the edge $z_1 v$  red,  a contradiction.  Thus $B=\{z_1, z_2\}$ and so $|C|\ge |B|=2$. 
Recall that   $e(G_r)\le 2n-3$.  By   \cref{l:K3-sat},    we have $ |C| \leq 3$  and $e(G_r) \ge 2n-4$.  \end{proof}

\noindent {\bf Claim\refstepcounter{counter}\label{C}  \arabic{counter}.} 
 $C = \{y_1,   y_2, y_3\}$.

\begin{proof}  Suppose  $|C|\ne 3$. By Claim~\ref{B}, we have $B=\{z_1, z_2\}$ and  $C = \{y_1,   y_2\}$.  Then   $|A|\ge 7$. Let $z_1'\in A$ be such that $z_1z_1'\notin E_b$. Then $z_1z_1'\notin E(G)$.      Suppose first   $yz\in E_b$. Since $\Delta(G_b)\le 2$, we may assume that $yy_1, zz_1\notin E_b$.  We may further assume that $d_b(y)\le d_b(z)$ in this case.    If  $yy_2\notin E_b$,  then we obtain a critical coloring of $G+z_1z_1'$ from $\tau$ by first coloring the edge $ z_1z_1'$ red, and then recoloring  $yz_1$ blue and all  edges   $z_1w$ red for  any $w\in N_b(z_1)\cap A$, a contradiction.  Thus   $yy_2\in E_b$, and so  $ zz_2\in E_b$ because $d_b(y)\le d_b(z)$. But then   we  obtain a critical coloring of $G+ z_1z_1'$ from $\tau$ by first coloring the edge $z_1z_1'$ red, and then recoloring  edges $yy_2$, $zz_2$ and all edges  incident with  $z_1$ or $ y_1$  in  $G_b$     red,     and  finally recoloring edges  $y_1z_1$, $y_2z_2$, $yz_1$, and $zy_1$ blue, a contradiction.    This proves that  $yz\notin E_b$.   
 Then $yz\notin E(G)$.  By  \cref{l:blue}(\ref{alpha}), either  $d_b(y)=2$ or  $d_b(z)=2$. Since $|B|=|C|=2$, we may assume that $d_b(y)=2$.  Then  $yy_1, yy_2\in E_b$. 
Suppose $d_b(z)\le 1$. We may assume that $zz_1\notin E_b$.  But then we obtain a critical coloring of $G+z_1z_1'$ from $\tau$  by first coloring the edge $z_1z_1'$ red, and then:   if $d_b(z)=0$ (and so $zz_2\notin E(G)$), recoloring $yz_1$ and $yz_2$ blue, and finally recoloring $yy_1$, $yy_2$, and all edges between $A$ and $B$ in $G_b$ red; if $d_b(z)=1$ (and so $zz_2\in E_b$),  recoloring edges $yz_1$, $yz_2$, $zy_1$, and $zy_2$ blue,  and  finally recoloring $yy_1$, $yy_2$, $zz_2$, and all edges between $z_1$ and $A$ in $G_b$ red.  This proves that  $d_b(z)=2$.  Then $zz_1, zz_2\in E_b$. Note that  $z_1z_2\in E_b$, else $z_1z_2\notin E(G)$, and we obtain a critical coloring of $G+ z_1z_1'$ from $\tau$ by first  coloring the edge $ z_1z_1'$ red, and then recoloring edges $yy_1$, $yy_2$, and all edges between $A$ and $B$ in $G_b$ red, and finally recoloring edges $yz_1$ and $yz_2$ blue.  By a similar argument, $y_1y_2\in E_b$. But then $G+yz$  has a critical coloring obtained from $\tau$ by first coloring   $yz$ blue, and then recoloring edges $z_1y$, $z_2y_1$, and $z_2y_2$ blue, and finally recoloring  $yy_1$, $yy_2$, $z_2z_1, z_2z$ red, a contradiction.\end{proof}

By Claim~\ref{B} and Claim~\ref{C},  $|B|=2$ and $|C|=3$.  Then $|A|\ge 6$.   By   \cref{l:K3-sat},  $e(G_r) = 2n-3$, and so  $e(G_b)=n-2$ by  Claim~\ref{Gb}.   Then $\delta(G_b)\le1$. Let $v_1, \ldots, v_s\in V(G)$ be all the vertices  satisfying $d_b(v_i)\leq 1$ for each $i\in [s]$. We may further assume that \[0=d_b(v_1)=\cdots= d_b(v_\ell)<1=d_b(v_{\ell+1})\le\cdots\le d_b(v_s),\] where $0\le \ell\le s$. Then $s-\ell$ is even. Note that   $(s-\ell)+2(n-s)=2e(G_b)=2n-4$, and so   $\ell+s=4$. Let $S:=\{v_1, \ldots, v_s\}$. By  \cref{l:blue}(\ref{alpha}),   $\alpha(G_b[S])\le2$.     Thus $\ell\le2$ and 
 \medskip

\noindent {\bf Claim\refstepcounter{counter}\label{vell}  \arabic{counter}.} 
       $\ell=s=2$, or $\ell=1$ and $s=3$, or $\ell = 0$ and $s=4$. 
\bigskip

\noindent {\bf Claim\refstepcounter{counter}\label{ScapA}  \arabic{counter}.} 
$S\cap A= \emptyset$.

\begin{proof}  Suppose $S\cap A\ne \emptyset$. Let $v\in A\cap S$. Then $d_b(v)\le1$. We may assume that $vz_1, vy_1 \notin E(G)$.  Note that if  $d_b(y)\le1$, then    we obtain a critical coloring of $G+vz_1$ from $\tau$ by first coloring the edge $vz_1$ red and then recoloring the edge $yv$ blue, a contradiction; Similarly, if $d_b(z)\le1$, then    we obtain a critical coloring of $G+vy_1$ from $\tau$ by first coloring the edge $vy_1$ red and then recoloring the edge $zv$ blue, a contradiction.
Thus  $d_b(y)=d_b(z)=2$. Then $S\subseteq A\cup B\cup C$. Recall that $d_b(v)\le 1$  and  $S$ is a clique in $G$. Since    both $A\cup B$ and $A\cup C $ are independent sets in $G_r$,      we see that either $B\cap S=\emptyset$ or $C\cap S=\emptyset$. It follows that either $S\subseteq A\cup B$ or $S\subseteq A\cup C$. In either case, $S$ is a clique  on at most two vertices in $G_b$,  and $\ell=0$ when $|S|=2$. Thus  $\ell+s\le 2$,    contrary to to the fact that $\ell+s=4$.     \end{proof}

\noindent {\bf Claim\refstepcounter{counter}\label{yz}  \arabic{counter}.} 
$y\in S$ or $z\in S$. 

\begin{proof} Suppose   $y,z\notin S$. Then $d_b(y)=d_b(z)=2$.  
 We may   assume that $yy_1, zz_1\in E_b$.  By Claim~\ref{ScapA},    $S\cap A= \emptyset$.  Then $S\subseteq  B\cup C$. Since  $zz_1\in E_b$ and  $  z_1z_2\notin E_r$, by  \cref{l:blue}(\ref{alpha}),        either $d_b(z_1)=2$ or $d_b(z_2)=2$.  
 Thus  $|B\cap S|\le1$.   We claim that $d_b(z_2)\ge 1$.  Suppose $d_b(z_2)=0$. 
 Then  $z_1z_2\notin E(G)$ and $d_b(z_1)=2$. Let $w, w^*\in A$ be  such that $z_1w\in E_b$ and $z_1w^*\notin E(G)$. But then we obtain a critical coloring of $G+z_1w^*$ from $\tau$ by first coloring the edge $z_1w^*$ red, and then recoloring edges $yz_1, z_2y_1$ blue and edges $yy_1, z_1w$ red, 
  a contradiction.  Thus $d_b(z_2)\ge 1$, as claimed. Then  $B \cap\{v_{ 1}, \ldots, v_\ell\} =\emptyset$.  Since $C$ is  an independent set in $G_r$, by  \cref{l:blue}(\ref{alpha}), $C\cap S$ is a clique on at most two vertices in $G_b$, and $\ell=0$ when $|S\cap C|=2$.   Thus  $\ell+s\le 3$,    contrary to   the fact that $\ell+s=4$.  
\end{proof}

By Claim~\ref{ScapA} and   Claim~\ref{yz}, we have   $S\cap A=\emptyset$  and $S\cap\{y,z\}\ne\emptyset$.   If $d_b(y)=d_b(z)=0$,   then $y,z\in S$ and $yz\notin E(G)$,  contrary to the fact that $S$ is a clique in $G$. Thus $\max\{d_b(y), d_b(z)\}\ge1$.
 Note that  $S\cap (B\cup C)\ne \emptyset$, else by Claim~\ref{vell}, $S=\{y,z\}$ and $d_b(y)=d_b(z)=0$, a contradiction.  
Suppose  $y\in S$. Then $d_b(y)\le1$. We may assume that $yy_2, yy_3\notin E_b$. Then $yy_2, yy_3\notin E(G)$. Thus $y_2, y_3\notin S$ because $S$ is a clique in $G$.    
  Then $B\cap S=\emptyset$, else, say $z_1\in S$, let $v\in A\setminus N_b(z_1)$. Then we obtain a critical coloring of $G+z_1v$ from $\tau$ by first coloring the edge $z_1v$ red and then recoloring the edge $yz_1$ blue, a contradiction.  Thus $\{y_1\}= S\cap (B\cup C)$ because $S\cap (B\cup C)\ne \emptyset$. Recall that $y\in S$ and $yy_1\notin E_r$. Since $S$ is a clique in $G$, we have   $yy_1\in E_b$. Then $d_b(y)=d_b(y_1)=1$ and $S\subseteq \{y, y_1, z\}$. By Claim~\ref{vell}, $s=3$, $\ell=1$ and  $z\in S$ with $d_b(z)=0$.  Then $yz\notin E(G)$, contrary to the fact that $S$ is a clique in $G$. This proves that $y\notin S$, and so   $z\in S$ by Claim~\ref{yz}. We may assume that $zz_2\notin E_b$. Then $zz_2\notin E(G)$ and so $z_2\notin S$. 
  Then $C\cap S=\emptyset$, else, say $y_1\in S$, let $v\in A\setminus N_b(y_1)$. Then we obtain a critical coloring of $G+y_1v$ from $\tau$ by first coloring the edge $y_1v$ red and then recoloring the edge $zy_1$ blue, a contradiction.  Thus $\{z_1\}= S\cap (B\cup C)$ because $S\cap (B\cup C)\ne \emptyset$. Then $S=\{z, z_1\}$. Since $S$ is a clique in $G$, we see that    $zz_1\in E_b$.  But then $\ell=0$ and $s=2$, contrary to   Claim~\ref{vell}. This completes the proof of  \cref{t:threeclaw}. \end{proof} 

\section*{Acknowledgements}
The authors would like to thank Jingmei Zhang for helpful discussion.

\end{document}